\documentclass[12pt]{article}

\usepackage{amsmath}
\usepackage{amssymb}
\usepackage{amsthm}

\usepackage{times}

\usepackage[T1]{fontenc}
\usepackage[polish]{babel}
\fontencoding{T1}
\usepackage[utf8]{inputenc}

\usepackage{amsfonts}
\usepackage{amsmath}
\usepackage{amssymb}
\usepackage{amsthm}

\usepackage{graphicx}
\usepackage{color}

\usepackage{pst-node,pst-plot,pst-eucl,pst-solides3d}

\author{Piotr Śniady (Wrocław)}
\date{}
\title{Asymptotyczna teoria reprezentacji grup permutacji}

\newcommand{\R}{\mathbb{R}}

\newcommand{\N}{\mathbb{N}}

\newcommand{\M}{{\mathcal{M}}}
\newcommand{\alert}[1]{\emph{#1}}

\DeclareMathOperator{\Tr}{Tr}
\DeclareMathOperator{\Id}{Id}
\DeclareMathOperator{\End}{End}

\newcommand{\Sym}[1]{\mathfrak{S}(#1)}
\newcommand{\Alt}[1]{\mathfrak{A}(#1)}

\theoremstyle{plain}
\newtheorem{lemma}{Lemat}
\newtheorem{theorem}[lemma]{Twierdzenie}

\newtheorem{wniosek}[lemma]{Wniosek}

\newtheorem{conjecture}[lemma]{Przypuszczenie}
\newtheorem{problem}[lemma]{Problem}

\begin{document}

\polishzx

\maketitle

\section{Wstęp}

Asymptotyczna teoria reprezentacji jest teorią stosunkowo młodą, gdyż została
ona zapoczątkowana pod koniec lat siedemdziesiątych XX wieku.  Choć formalnie
jest ona częścią ,,zwykłej'' teorii reprezentacji, to wyróżnia się ona paletą
używanych środków. Asymptotyczna teoria reprezentacji oprócz zwykłych
metod teorii reprezentacji wykorzystuje metody analizy i nowe metody
kombinatoryczne. Ponadto blisko jest związana z tak ostatnio modnym działem
matematyki, jakim jest teoria macierzy losowych oraz z \emph{wolną
probabilistyką} Dana Voiculescu. W niniejszym artykule chciałbym przedstawić
spojrzenie z lotu ptaka na tę szybko rozwijającą się teorię.

Niniejszy tekst jest znacznie rozszerzonym zapisem XV Wykładu im.~Wojtka
Pulikowskiego, który wygłosiłem 30 maja 2008 r.~na Wydziale Matematyki i
Informatyki Uniwersytetu Adama Mickiewicza w Poznaniu. Przewidując obecność
licznej grupy poznańskich kombinatoryków podczas wykładu, zdecydowałem się w
szczególny sposób wyeksponować właśnie kombinatoryczne zagadnienia i problemy.
Moim celem było pokazanie, że klasyczne metody
kombinatoryczne często zawodzą w problemach asymptotycznej teorii reprezentacji,
a teoria ta tworzy nowe kombinatoryczne narzędzia. W
niniejszym artykule pozwolę sobie przedstawić temat nieco szerzej i
zaprezentuję również kontekst zarówno klasycznej jak i asymptotycznej teorii
reprezentacji.

\section{Teoria reprezentacji}

\subsection{Grupy}
Niezliczone są sytuacje w matematyce i fizyce, gdy rozwiązanie problemu staje
się znacznie prostsze dzięki uwzględnieniu jego \emph{symetrii}. Pojęcie zbioru
symetrii danego obiektu okazało się tak ważne, że zostało ono sformalizowane
pod nazwą \emph{grupy}, a \emph{teoria grup} stała się jedną z najważniejszych
działów matematyki. 

\begin{figure}
\begin{center}
\psset{unit=10ex}
\begin{pspicture}(-1,-0.33)(1,1)
\psaxes[labels=none,ticks=no,linecolor=gray]{->}(0,0)(-1,-0.5)(1,1)
\cnode*[linecolor=gray](0,0.577350269){1mm}{a}
\cnode*[linecolor=gray](0.5,-0.288675135){1mm}{b}
\cnode*[linecolor=gray](-0.5,-0.288675135){1mm}{c}
\ncline[linecolor=gray]{-}{a}{b}
\ncline[linecolor=gray]{-}{b}{c}
\ncline[linecolor=gray]{-}{c}{a}
\rput(0.2,0.57735){$1$}
\rput(0.65,-0.3){$2$}
\rput(-0.65,-0.3){$3$}
\end{pspicture}
\end{center}
\caption{Trójkąt równoboczny na płaszczyźnie, którego środek ciężkości znajduje
się w początku układu współrzędnych.}
\label{rys:trojkat}
\end{figure}
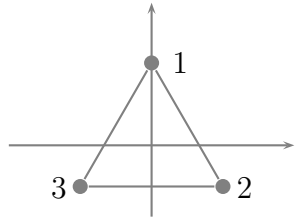

Początkowo obiektem zainteresowania  matematyków były konkretne grupy rozumiane
właśnie jako zbiory symetrii bardzo konkretnych obiektów o charakterze
geometrycznym, kombinatorycznym lub algebraicznym. Przykładem takiego podejścia
jest badanie symetrii trójkąta równobocznego widocznego na Rysunku
\ref{rys:trojkat}. W tym prostym przykładzie, który jeszcze będzie przedmiotem
dalszej analizy, jeśli układ współrzędnych dobierzemy w ten sposób, aby jego
początek znajdował się w środku masy trójkąta, wówczas każda izometria
rozważanego trójkąta jest \emph{odwzorowaniem liniowym}, czyli daje się opisać
przy pomocy pewnej macierzy, a składanie izometrii odpowiada mnożeniu
wspomnianych macierzy. Innymi słowy, grupą symetrii naszego trójkąta jest
pewien zbiór odwzorowań liniowych (macierzy).

Z biegiem czasu obiektem zainteresowania matematyków stały się grupy widziane w
sposób coraz bardziej abstrakcyjny. O grupie coraz mniej myślano jak o
konkretnym zbiorze symetrii jakiegoś obiektu, a coraz bardziej jak o
abstrakcyjnym zbiorze spełniającym pewne aksjomaty. Aby dać posmak takiego
bardziej abstrakcyjnego myślenia oznaczmy wierzchołki
trójkąta rozważanego w naszym przykładzie liczbami $1$, $2$, $3$. Wówczas
wspomniana grupa symetrii trójkąta
równobocznego może być widziana jako grupa permutacji wierzchołków
trójkąta czyli zbioru $\{1,2,3\}$. Z tego
powodu grupę tę będziemy oznaczać symbolem $\Sym{3}$ (dalszą dyskusję tej 
notacji
przedstawimy nieco dalej).

\subsection{Reprezentacje}
Mówiąc obrazowo, przedmiotem \emph{teorii reprezentacji} jest badanie
sposobów, w jakie abstrakcyjne grupy realizują się w konkretny sposób jako
symetrie prostych obiektów geometrycznych, a więc w pewnym sensie stanowi ona
powrót abstrakcyjnej teorii grup do jej konkretnych korzeni. Patrząc na
przykład grupy $\Sym{3}$ z punktu widzenia teorii reprezentacji moglibyśmy
powiedzieć (być może nadużywając terminologii), że grupa $\Sym{3}$
reprezentuje się jako izometrie trójkąta.

Formalna definicja brzmi następująco: \alert{reprezentacją} grupy $G$ nazywamy
odwzorowanie $\rho$, które przyporządkowuje elementom grupy odwracalne macierze
(o jakimś ustalonym rozmiarze $n$)
$$\rho:G\rightarrow \M_{n}$$
i które jest \emph{homomorfizmem}, to znaczy iloczynowi elementów grupy
powinien
odpowiadać iloczyn odpowiadających im macierzy:
$$ \rho(g_1 g_2)= \rho(g_1) \rho(g_2)\qquad \text{dla dowolnych }
g_1,g_2\in G.$$
Macierze, którymi się będziemy zajmować w niniejszym artykule, będą zawsze
kwadratowymi macierzami o wyrazach rzeczywistych lub zespolonych.

Alternatywnie, możemy myśleć, że reprezentacja przyporządkowuje elementom grupy
odwracalne odwzorowania liniowe pewnej ustalonej skończenie wymiarowej
przestrzeni liniowej $V$:
$$\rho:G\rightarrow \End(V)$$
 oraz że
iloczynowi elementów grupy ma odpowiadać złożenie odpowiednich przekształceń:
$$ \rho(g_1 g_2)= \rho(g_1)\circ \rho(g_2)\qquad \text{dla dowolnych }
g_1,g_2\in G.$$

\subsection{Przykład: reprezentacje grupy $\Sym{3}$}
\label{subsec:s3}

Z jedną reprezentacją grupy $\Sym{3}$ już się spotkaliśmy: jest to
reprezentacja,
która permutacjom zbioru $\{1,2,3\}$ przyporządkowuje odpowiednią izometrię
trójkąta.

Są też inne reprezentacje tej
grupy: możemy na przykład każdemu elementowi $\Sym{3}$ przypisać odwzorowanie
identycznościowe prostej (innymi słowy, jest to macierz $[1]$ o rozmiarach
$1\times 1$); Czytelnik łatwo domyśli się, dlaczego reprezentacja ta nazywana
jest \emph{reprezentacją trywialną}. 

Jeszcze inną reprezentacją jest \emph{reprezentacja alternująca}, która
permutacjom parzystym przyporządkowuje odwzorowanie identycznościowe na prostej,
czyli odwzorowanie o macierzy $[1]$, a permutacjom nieparzystym
przyporządkowuje symetrię prostej względem $0$, czyli odwzorowanie o macierzy
$[-1]$.

\subsection{Przykład: reprezentacja grupy $\Alt{5}$}
\label{subsec:a5}

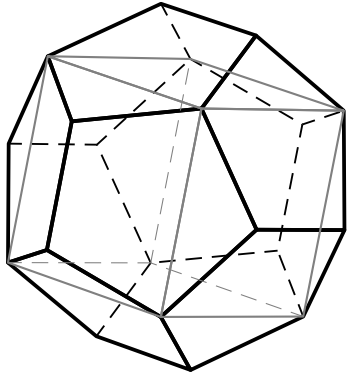
\begin{figure}
\psset{unit=3.5cm}
\psset{SphericalCoor,viewpoint=100 13 -25,Decran=50}
\begin{center}
\begin{pspicture}(-0.7,-0.7)(0.7,0.7)
{\psSolid[linewidth=0.5mm,action=draw,object=new,sommets=
0 0 0
-0.5 0.0 -1.30902
-0.5 0.0 1.30902
0.0 -1.30902 -0.5
0.0 -1.30902 0.5
0.0 1.30902 -0.5
0.0 1.30902 0.5
0.5 0.0 -1.30902
0.5 0.0 1.30902
-1.30902 -0.5 0.0
-1.30902 0.5 0.0
-0.809017 -0.809017 -0.809017
-0.809017 -0.809017 0.809017
-0.809017 0.809017 -0.809017
-0.809017 0.809017 0.809017
0.809017 -0.809017 -0.809017
0.809017 -0.809017 0.809017
0.809017 0.809017 -0.809017
0.809017 0.809017 0.809017
1.30902 -0.5 0.0
1.30902 0.5 0.0,
faces={
[14 10 9 12 2]
[7 17 20 19 15]
[15 19 16 4 3]
[19 20 18 8 16]
[20 17 5 6 18]
[17 7 1 13 5]
[7 15 3 11 1]
[6 5 13 10 14]
[13 1 11 9 10]
[11 3 4 12 9]
[4 16 8 2 12]
[8 18 6 14 2]}
]}
{\psSolid[linecolor=gray,action=draw,object=new,sommets=
0 0 0
-0.5 0.0 -1.30902
-0.5 0.0 1.30902
0.0 -1.30902 -0.5
0.0 -1.30902 0.5
0.0 1.30902 -0.5
0.0 1.30902 0.5
0.5 0.0 -1.30902
0.5 0.0 1.30902
-1.30902 -0.5 0.0
-1.30902 0.5 0.0
-0.809017 -0.809017 -0.809017
-0.809017 -0.809017 0.809017
-0.809017 0.809017 -0.809017
-0.809017 0.809017 0.809017
0.809017 -0.809017 -0.809017
0.809017 -0.809017 0.809017
0.809017 0.809017 -0.809017
0.809017 0.809017 0.809017
1.30902 -0.5 0.0
1.30902 0.5 0.0,
faces={
[2 6 13 9]
[2 9 3 16]
[2 16 20 6]
[3 7 20 16]
[3 9 13 7]
[6 20 7 13]}
]}
\end{pspicture}
\end{center}
\caption{Dwunastościan foremny z sześcianem wpisanym w taki sposób, że każdy
wierzchołek sześcianu jest również wierzchołkiem dwunastościanu.}
\label{rys:dwunastoscian}
\end{figure}

Powyższy przykład jest tak prosty, że mógł wywołać u Czytelnika błędne wrażenie,
jakoby reprezentacje były czymś trywialnym. Aby zatrzeć to wrażenie przedstawię
teraz znacznie mniej  oczywisty przykład. 

W dwunastościan foremny można wpisać sześcian w ten sposób aby każdy wierzchołek
sześcianu był również wierzchołkiem dwunastościanu. Jeden ze sposobów wpisania
takiego sześcianu przedstawiony jest na Rysunku~\ref{rys:dwunastoscian}. W sumie
takich sześcianów jest pięć.

Wynika stąd, że każda izometria dwunastościanu zachowująca
orientację przestrzeni (innymi słowy: każdy obrót dwunastościanu wokół pewnej
osi przekształcający dwunastościan w siebie) zadaje pewną permutację powyższych
pięciu sześcianów. Można udowodnić, że powyższa permutacja pięciu sześcianów
jest zawsze permutacją parzystą (czyli daje się zapisać jako iloczyn parzystej
liczby transpozycji). Permutacje parzyste zbioru pięcioelementowego tworzą grupę
oznaczaną symbolem $\Alt{5}$. 

Co więcej, można wykazać, że
powyższa odpowiedniość jest wzajemnie jednoznaczna, a zatem każdemu elementowi
$\Alt{5}$ odpowiada dokładnie jedna izometria dwunastościanu zachowująca
orientację.
Powinno być w miarę jasne, że odwzorowanie to jest homomorfizmem, czyli
iloczynowi dwóch permutacji odpowiada złożenie odpowiadających im izometrii.

Jeśli układ współrzędnych wybierzemy tak, aby jego początek znajdował się w
środku masy dwunastościanu, to wspomniane izometrie są odwzorowaniami liniowymi,
a zatem odwzorowanie, które parzystym permutacjom pięciu sześcianów
przyporządkowuje izometrię liniową dwunastościanu, jest reprezentacją grupy
$\Alt{5}$ na przestrzeni trójwymiarowej.

\subsection{Reprezentacje nieredukowalne}

Jeśli $\rho_1:G\rightarrow \M_{n_1}$ oraz $\rho_2:G\rightarrow \M_{n_2}$
są
reprezentacjami tej samej grupy $G$, możemy zdefiniować nową reprezentację
grupy $G$, zwaną \emph{sumę prostą} $\rho_1\oplus\rho_2:G\rightarrow
\M_{n_1+n_2}$,
która jest zadana przez macierze blokowe.
$$ (\rho_1\oplus\rho_2)(g) = \begin{bmatrix} \rho_1(g) &  \\ &\rho_2(g)
\end{bmatrix}. $$
Powyższą definicję możemy też sformułować następująco: jeśli
$\rho_1:G\rightarrow \End(V_1)$ oraz $\rho_2:G\rightarrow \End(V_2)$ są
reprezentacjami tej samej grupy $G$ na przestrzeniach liniowych $V_1$ i $V_2$,
to $\rho_1\oplus\rho_2:G\rightarrow \End(V_1\oplus V_2)$ jest reprezentacją
zadaną wzorem
$$ (\rho_1\oplus\rho_2)(g)= \big(\rho_1(g),\rho_2(g)\big). $$

Reprezentacje będące sumami prostymi są zwykle mniej interesujące od tych,
które nie są tej postaci. Możemy więc zapytać, czy zadana reprezentacja
nie jest sumą prostą mniejszych
reprezentacji. To pojęcie sformalizowane jest w następujący sposób: mówimy, że
reprezentacja $\rho:G\rightarrow \End(V)$ na przestrzeni liniowej $V$ jest
\alert{redukowalna}, jeśli istnieje rozkład na podreprezentacje:
$V=V_1\oplus V_2$ oraz $\rho=\rho_1\oplus \rho_2$.  Dla pełnej ścisłości
powinniśmy jeszcze wymagać, by przestrzenie $V_1$ i $V_2$ były nietrywialne,
czyli nie składały się tylko z wektora zerowego. 

Mówimy też, że reprezentacja $\rho:G\rightarrow \End(V)$ jest
\emph{przywiedlna}, jeśli istnieje
podprzestrzeń $V_1$ która jest \emph{podprzestrzenią niezmienniczą}, to znaczy
dla dowolnego elementu $g\in G$, jeśli $v\in V_1$ to również $\rho(g) v\in V_1$.
Ponownie powinniśmy założyć, że $V_1$ jest nietrywialna, to znaczy że nie składa
się tylko z wektora zerowego ani nie jest równa całej przestrzeni $V$. W
interesującym nas w tym artykule przypadku grup skończonych powyższe dwie
własności: redukowalność i przywiedlność są równoważne, nie musimy ich więc
specjalnie rozróżniać.

Reprezentację, która nie jest redukowalna, nazywamy \emph{nieredukowalną}.
Reprezentacje nieredukowalne pełnią w teorii reprezentacji podobną rolę jak
liczby pierwsze w teorii liczb czy atomy w chemii, a zatem są elementarnymi
cegiełkami, z których zbudowane są wszystkie reprezentacje, gdyż
każdą reprezentację można rozłożyć jako sumę prostą pewnej liczby reprezentacji
nieredukowalnych (ponieważ w niniejszym artykule interesują nas tylko grupy
skończone i wyłącznie reprezentacje na skończenie wymiarowych przestrzeniach
liniowych, możemy nie martwić się pewnymi patologicznymi sytuacjami).

Jak się okazuje, wszystkie przykłady przedstawione w Rozdziałach \ref{subsec:s3}
i \ref{subsec:a5} są reprezentacjami nieredukowalnymi. Z kolei reprezentacja
grupy $\Sym{3}$ na przestrzeni $\R^3$, w której
permutacje działają na wektorach przez permutowanie ich współrzędnych, jest
rozkładalna, gdyż jest sumą reprezentacji trywialnej oraz dwuwymiarowej
reprezentacji z Rozdziału \ref{subsec:s3}.

Grupa $\Sym{3}$ okazuje się nie mieć innych reprezentacji
nieredukowalnych niż te przedstawione w Rozdziale \ref{subsec:s3}. Zdanie to
wymaga pewnego doprecyzowania, gdyż poprzez wybór innego układu współrzędnych
dla trójkąta na płaszczyźnie na Rysunku~\ref{rys:trojkat} możemy uzyskać
nieskończenie wiele reprezentacji nieredukowalnych grupy $\Sym{3}$. Otóż
umawiamy
się, że dwie reprezentacje tej samej grupy są \emph{równoważne} (takich
reprezentacji nie będziemy w przyszłości rozróżniać), jeśli poprzez zmianę
układów współrzędnych w odpowiednich przestrzeniach liniowych odpowiednie
macierze stają się sobie równe.

\subsection{Zastosowanie teorii reprezentacji w fizyce i chemii}

\begin{figure}
\begin{center}
\psset{unit=7ex}
\begin{pspicture}(-2,-1.6)(2,1.6)
\SpecialCoor
\degrees[360]
\cnodeput*(1;0){c1}{C}
\cnodeput*(1;60){c2}{C}
\cnodeput*(1;120){c3}{C}
\cnodeput*(1;180){c4}{C}
\cnodeput*(1;240){c5}{C}
\cnodeput*(1;300){c6}{C}
\cnodeput*(1.6;0){d1}{H}
\cnodeput*(1.6;60){d2}{H}
\cnodeput*(1.6;120){d3}{H}
\cnodeput*(1.6;180){d4}{H}
\cnodeput*(1.6;240){d5}{H}
\cnodeput*(1.6;300){d6}{H}
\ncline{c1}{c2}
\ncline{c2}{c3}
\ncline{c3}{c4}
\ncline{c4}{c5}
\ncline{c5}{c6}
\ncline{c6}{c1}
\ncline{c1}{d1}
\ncline{c2}{d2}
\ncline{c3}{d3}
\ncline{c4}{d4}
\ncline{c5}{d5}
\ncline{c6}{d6}
\pscircle(0,0){0.6}
\end{pspicture}
\end{center}
\caption{Schemat cząsteczki benzenu.}
\label{rys:benzene}
\end{figure}
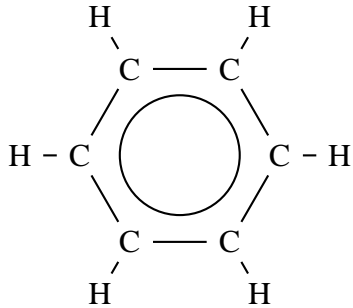

Teoria reprezentacji jest bardzo wygodnym narzędziem dla fizyków oraz chemików
badających układy kwantowe o dużej symetrii. Do opisu takiego układu używamy
\emph{przestrzeni Hilberta} $\mathcal{H}$, która jest pewną nieskończenie
wymiarową przestrzenią liniową.  Zakładamy, że nasz układ kwantowy ma
nietrywialną grupę symetrii $G$, na przykład dla cząsteczki
benzenu (Rysunek \ref{rys:benzene}) grupą $G$ jest grupa skończona symetrii
sześciokąta foremnego.

Oddziaływania w naszym układzie kwantowym opisywane są przez \emph{hamiltonian}
$H:\mathcal{H}\rightarrow\mathcal{H}$, który jest pewnym odwzorowaniem liniowym
na przestrzeni Hilberta $\mathcal{H}$ czy też, mówiąc nieco innym językiem,
jest pewną nieskończoną macierzą. W fizyce i chemii kwantowej częste jest
pytanie o dopuszczalne poziomy energii interesującego nas układu kwantowego. W
języku matematycznym problem ten tłumaczy się jako wyznaczenie wartości własnych
hamiltonianu. 

Wyznaczenie wspomnianych wartości własnych nie jest prostym problemem, nawet
jeśli chcemy się zadowolić przybliżoną odpowiedzią uzyskaną dzięki
dyskretyzacji i zastąpieniu nieskończenie wymiarowej przestrzeni Hilberta przez
przestrzeń o skończonym wymiarze i zastosowaniu metod komputerowych. Przyczyna
tej trudności jest następująca: wymiar przestrzeni Hilberta rośnie wykładniczo
wraz z liczbą analizowanych cząstek i uwzględnienie nawet niewielkiej liczby
elektronów powoduje, że musimy manipulować bardzo dużymi macierzami.

Zauważmy jednak, że każdemu elementowi grupy symetrii $G$ naszej cząsteczki
odpowiada pewne przekształcenie przestrzeni Hilberta $\mathcal{H}$, które
opisuje w jaki sposób dane geometryczne przekształcenie realizuje się w układzie
fizycznym. Innymi słowy, otrzymaliśmy reprezentację grupy $G$ na
przestrzeni $\mathcal{H}$. Nasze założenie, że grupa $G$ jest grupą symetrii
cząsteczki oznacza, że również hamiltonian $H$ opisujący dynamikę cząsteczki
jest niezmienniczy na działanie grupy $G$, a zatem jego przestrzenie własne są
podprzestrzeniami niezmienniczymi dla naszej reprezentacji. W szczególności
każda przestrzeń własna hamiltonianu może być traktowana jako reprezentacja
grupy $G$. Wynika stąd, że zamiast badać hamiltonian jako odwzorowanie liniowe
na dużej przestrzeni Hilberta, możemy dla każdej reprezentacji nieredukowalnej
$\rho$ grupy $G$ wybrać tę część przestrzeni Hilberta $\mathcal H$, na której
grupa $G$ reprezentuje się w sposób analogiczny do reprezentacji $\rho$ i
dalej badać hamiltonian na tej mniejszej przestrzeni. Dzięki temu musimy
wyznaczyć wartości własne znacznie mniejszych macierzy niż bez uwzględnienia
symetrii problemu.

W przypadku układu kwantowego, który może swobodnie rotować w przestrzeni wokół
swojego środka masy, grupą symetrii jest grupa obrotów. Jak się okazuje,
reprezentacje nieredukowalne tej grupy
odpowiadają stanom kwantowym o ustalonej wartości momentu pędu. Czysto
matematyczna analiza tych reprezentacji pozwala uzyskać na przykład formuły na
dodawanie kwantowego momentu pędu w piękny, abstrakcyjny sposób bez
jakiegokolwiek odwoływania się do szczegółów interesujących nas układów
kwantowych.

\subsection{Charaktery}

Jeśli o reprezentacji $\rho:G\rightarrow \M_n$ myślimy jak o funkcji, która
elementom grupy przyporządkowuje macierze, to zmiana układu współrzędnych może
spowodować zmianę wyrazów tej macierzy. Skoro w dwóch układach współrzędnych ta
sama reprezentacja może wyrażać się przez inne macierze, oznacza to, że być
może dobrym pomysłem byłoby znalezienie innych wielkości opisujących
reprezentację, a które byłyby niezależne od wyboru układu współrzędnych. 

Przykładem takiej niezmienniczej wielkości jest \emph{charakter reprezentacji}
zdefiniowany wzorem 
$$ \chi^\rho(g) = \Tr \rho(g),$$
gdzie $\Tr$ oznacza ślad macierzy. To nieco zaskakujące, ale jak się okazuje,
niemal wszystkie pytania teorii reprezentacji dają się przeformułować na
pytania dotyczące charakterów reprezentacji.

\subsection{Zastosowanie charakterów: ile razy trzeba przetasować talię kart?}
\label{subsec:tasowanie}

Jako zastosowanie teorii reprezentacji przeanalizujemy następujący problem: ile
razy należy potasować talię kart, aby była ona dobrze przetasowana? Nieco
bardziej formalnie: niech zmiana pozycji kart w kolejnych tasowaniach będzie
opisana przez ciąg niezależnych zmiennych losowych o jednakowym rozkładzie
$X_1,X_2,\dots$. Wspomniane zmienne losowe przyjmują wartości w grupie
permutacji $\Sym{n}$, gdzie w większości gier karcianych $n=52$. Interesuje nas
rozkład zmiennej losowej $X_1\circ \cdots \circ X_k$ opisującej rozkład kart po
$k$ przetasowaniach oraz jak odległy jest ten rozkład od rozkładu
jednostajnego na grupie permutacji $\Sym{n}$, który odpowiada idealnemu
przetasowaniu kart.

Niech $\mu$ będzie miarą probabilistyczną na $\Sym{n}$, która opisuje rozkład
zmiennych losowych $X_1,X_2,\dots$. O $\mu$ możemy myśleć jak o funkcji, która
elementom grupy przypisuje ich prawdopodobieństwa, ale też jak o elemencie
algebry grupowej $\R[\Sym{n}]$, czyli jak o formalnej kombinacji liniowej
elementów grupy (współczynnik stojący przy permutacji jest równy
prawdopodobieństwu tejże permutacji). Jeśli o rozkładzie zmiennej losowej
$X_1\circ\cdots\circ X_k$ myślimy jak o rozkładzie prawdopodobieństwa, to
jest on równy $\mu^{\ast k}$, czyli $k$-krotnemu splotowi $\mu$ ze
sobą; jeśli zaś myślimy o nim jak o elemencie algebry grupowej, to jest on
równy $\mu^k$, czyli $k$-krotnej potędze elementu $\mu\in\R[\Sym{n}]$. Oznacza
to,
że nasz problem sprowadza się do efektywnego obliczania potęg w algebrze
grupowej interesującej nas grupy $\Sym{n}$.

Jeśli $\rho$ jest reprezentacją grupy, wówczas
$$ \rho( \mu^k ) = \big[ \rho(\mu) \big]^k; $$
innymi słowy dowolna reprezentacja zamienia splot na grupie na
iloczyn macierzy. Nasza strategia polega na tym, aby użyć pewnej rodziny
$(\rho_\alpha)$ reprezentacji; jeśli nasza rodzina będzie dostatecznie duża,
mamy nadzieję zrekonstruować informację na temat $\mu^k$ z informacji na temat
$\rho_\alpha(\mu^k )$.

Gdyby zmienne losowe $X_1,\dots$ przyjmowały wartości liczbowe, a interesującą
nas wielkością byłaby suma $X_1+\cdots+X_k$, mielibyśmy do czynienia z
przemienną grupą $\R$ wyposażoną w dodawanie. Dla dowolnej liczby $z$
odwzorowanie
$$ \rho_z :  x \mapsto e^{ixz} $$
może być traktowane jako reprezentacja grupy $\R$ (o wartościach w macierzach o
rozmiarach $1\times 1$). Widać, że (pomijając pewne trudności związane z
przejściem od sytuacji dyskretnej do ciągłej) $\rho_z(\mu)$ jest transformatą
Fouriera miary probabilistycznej $\mu$, a zatem nasz program
jest uogólnieniem transformaty Fouriera na przypadek grup
nieprzemiennych.

Jak się okazuje, najlepszym wyborem rodziny reprezentacji $(\rho_\alpha)$ są
wszystkie reprezentacje nieredukowalne interesującej nas grupy. Z jednej strony
dzięki temu potęgować będziemy możliwie proste macierze, z drugiej strony
reprezentacji tych jest wystarczająco wiele, aby pozwolić na rekonstrukcję
informacji na temat interesującej nas miary $\mu^k$.

Szczególnie interesujący jest przypadek, kiedy miara $\mu$ jest
\emph{centralna}, to znaczy przypisuje stałe prawdopodobieństwo wszystkim
elementom grupy z tej samej klasy sprzężoności. Na przykład, jeśli losowo
wybierzemy dwie karty z talii i zamienimy je miejscami, rozkład tak uzyskanej
permutacji jest centralny. Jak się okazuje, w takim przypadku dla dowolnej
reprezentacji nieredukowalnej $\rho$ macierz $\rho(\mu)$ jest
skalarną wielokrotnością macierzy jednostkowej. Wielkość tego skalara można
wyznaczyć z równości
$$ \rho(\mu) = \frac{\Tr \rho(\mu) }{ \Tr \Id } \Id = \frac{\Tr \rho(\mu) }{ \Tr
\rho(e) } \Id, $$
gdzie wykorzystaliśmy fakt, że $e$, element neutralny grupy, reprezentuje się
jako macierz jednostkowa. Oznacza to, że w przypadku, gdy $\mu$ jest miarą
centralną na interesującej nas grupie, nie musimy się martwić potęgowaniem
macierzy, a jedynie znacznie prostszym potęgowaniem liczb
\begin{equation} 
\label{eq:normalizowane}
\rho \mapsto \frac{\Tr \rho(\mu) }{ \Tr
\rho(e) } = \frac{\chi^\rho(\mu)}{\chi^\rho(e)}, 
\end{equation}
podobnie jak w przypadku zwykłej
transformaty Fouriera.

Widać więc, że niekomutatywna transformata Fouriera, jakiej dostarcza teoria
reprezentacji do badania nieprzemiennych grup, zdefiniowana jest przy pomocy
charakterów i że do uprawiana analizy harmonicznej na grupach konieczne
jest dobre zrozumienie charakterów.


\subsection{Dalsza lektura}

Niniejszy przeglądowy artykuł nie ma
ambicji dostarczyć pełnej informacji bibliograficznej, dlatego prawie nie ma w
nim cytowań. Czytelnika zainteresowanego dalszą lekturą zachęcam do książki
Serre'a \cite{SerreReprezentacje} (dostępnej także w polskim tłumaczeniu)
będącej doskonałym wstępem do teorii reprezentacji grup
skończonych. O zastosowaniu teorii reprezentacji do badania tasowania kart
stanowi artykuł Diaconisa i Shahshahaniego \cite{DiaconisShahshahani1981} oraz 
monografia Diaconisa \cite{Diaconis1988}.

\section{Asymptotyczna teoria reprezentacji grup permutacji: charaktery}

\subsection{Asymptotyczna teoria reprezentacji}

Mówiąc w wielkim skrócie, jak przedmiotem badań teorii reprezentacji jest grupa
$G$ i jej reprezentacja $\rho$, tak przedmiotem badań asymptotycznej teorii
reprezentacji jest ciąg grup $G_1,G_2,\dots$ oraz odpowiadający im ciąg
reprezentacji $\rho_1,\rho_2,\dots$ (to znaczy: $\rho_n$ jest reprezentacją
grupy $G_n$). W asymptotycznej teorii reprezentacji pytamy: \emph{co możemy
powiedzieć o reprezentacji $\rho_n$ w granicy, gdy $n\to\infty$?}

Aby na tak sformułowane pytanie dało się cokolwiek odpowiedzieć, reprezentacje
grup $G_1,G_2,\dots$ z rozważanego ciągu muszą mieć pewną wspólną strukturę,
pozwalającą w jakiś sposób porównywać ze sobą reprezentacje różnych grup. Z
tego powodu w dalszym ciągu rozważać będę tylko ciąg grup permutacji
$\Sym{1},\Sym{2},\dots$, których reprezentacje taką wspólną strukturę mają.
Ponadto, jest to sam w sobie interesujący ciąg grup, gdyż każda grupa skończona
jest podgrupą dostatecznie dużej grupy permutacji. Co więcej, grupy alternujące
$\Alt{n}$ (złożone z parzystych permutacji, a zatem bardzo blisko związane z
grupami $\Sym{n}$) są \emph{grupami prostymi} dla $n\geq 5$, przez co są
szczególnie ważne jako przykłady elementarnych cegiełek, z których zbudowane są
wszystkie grupy skończone.

\subsection{Nieredukowalne reprezentacje grup permutacji $\Sym{n}$}

\begin{figure}[tb]
\begin{center}
\includegraphics
{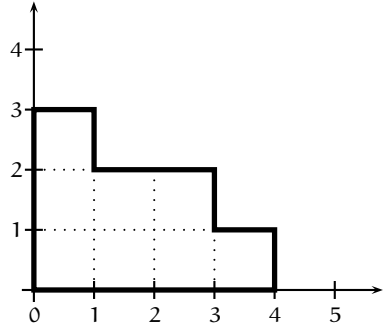}
\end{center}
\caption{Diagram Younga $(4,3,1)$ narysowany w konwencji francuskiej.}
\label{rys:francuski}
\end{figure}

Jak się okazuje, nieredukowalne reprezentacje grupy $\Sym{n}$ są we wzajemnie
jednoznacznej odpowiedniości z \emph{diagramami Younga} o $n$ klatkach. 
Przykład takiego diagramu widoczny jest na Rysunku \ref{rys:francuski}. Bez
użycia
środków graficznych diagram Younga $\lambda$ może być zdefiniowany jako
nierosnący ciąg $\lambda=(\lambda_1\geq \cdots \geq \lambda_l)$, którego
wyrazami są liczby całkowite dodatnie. Alternatywnie, dopisując do takiego ciągu
nieskończenie wiele zer, diagram Younga możemy zdefiniować jako nierosnący ciąg 
ciąg $\lambda=(\lambda_1\geq \lambda_2\geq \cdots)$, którego
wyrazami są liczby całkowite nieujemne i który składa się tylko ze
skończenie wielu niezerowych elementów. 

Element $\lambda_i$
interpretujemy jako liczbę klatek w $i$-tym wierszu, zatem liczba klatek
diagramu oznaczana przez $|\lambda|$ spełnia
$|\lambda|=\lambda_1+\lambda_2+\cdots$. Terminy: \emph{liczba wierszy} oraz
\emph{liczba kolumn} diagramu Younga $\lambda$ powinny być jasne w podejściu
graficznym do diagramów Younga; alternatywnie można je zdefiniować, odpowiednio,
jako $\lambda_1$ oraz jako liczbę niezerowych wyrazów ciągu
$\lambda_1,\lambda_2,\dots$.

Niestety, szczegóły konstrukcji nieredukowalnej reprezentacji $\rho^\lambda$
odpowiadającej diagramowi $\lambda$ nie są proste. Klasyczna metoda ich
konstrukcji oparta na \emph{symetryzatorach Younga} jest mało intuicyjna, a
diagramy Youga pojawiają się w niej \emph{deus ex machina}. Na szczęście w
ostatnich latach dostępna stała się piękna i naturalna metoda oparta na
\emph{elementach Jucysa-Murphyego}; zainteresowanego
Czytelnika odsyłam do doskonałej pracy Okounkova i Vershika
\cite{OkounkovVershik1996}.

\subsection{Charaktery grup permutacji i reguła Murnaghana-Na\-ka\-ya\-my}

Co prawda skonstruowanie nieredukowalnych reprezentacji grup permutacji
jest niełatwe, jednak metoda obliczania charakterów---które są
niemal jedyną po\-trze\-bną  nam do zastosowań wielkością---jest stosunkowo
prosta, zadana przez \emph{regułę Mur\-na\-gha\-na-Nakayamy}. Poniżej
przedstawię tę regułę w jednym tylko celu: aby pokazać, jak bardzo nie nadaje
się ona do pewnych problemów asymptotycznej teorii reprezentacji (gdyż do
pewnych rodzajów oszacowań asymptotycznych nadaje się ona bardzo dobrze), tak
więc Czytelnik bez straty ciągłości może opuścić niniejszy rozdział i
kontynuować lekturę w rozdziale \ref{sec:nastepny-po-murnaghan}.

Reguła Murnaghana-Nakayamy głosi,
że aby obliczyć wartość charakteru $\chi^\lambda(\pi)= \Tr \rho^\lambda(\pi)$,
powinniśmy najpierw wyznaczyć rozkład na cykle permutacji $\pi$ oraz wyznaczyć
długości $l_1,\cdots,l_k$ jej cykli. Dla przykładu,
permutacja $\pi=(3,5)(2,4)(7,1,6,8)$ składa się z trzech cykli o długościach
odpowiednio $2,2,4$ (kolejność, w jakiej ustawiliśmy cykle, nie jest istotna).

\begin{figure}[tb]
\begin{center}
\includegraphics[width=0.3\textwidth]{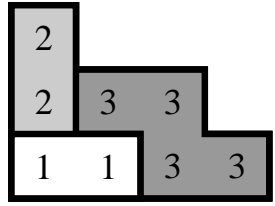}
\end{center}
\caption{Przykład podziału diagramu $(4,3,1)$ z Rys.~\ref{rys:francuski} na
skośne
paski o długościach: $2,2,4$. Kolejne paski zostały ponumerowane i oznaczone
kolorami począwszy od białego ku szarości.}
\label{rys:murnaghan}
\end{figure}

Następnie powinniśmy znaleźć wszystkie rozkłady diagramu $\lambda$ na
\emph{skośne paski} o długościach $l_1,\dots,l_k$. Przykład takiego rozkładu
widoczny jest na Rysunku~\ref{rys:murnaghan}. Każdy skośny pasek to kolekcja
klatek
diagramu Younga o tej własności, że można przejść wszystkie klatki paska
poruszając się tylko w prawo lub w dół o jedno pole. 
Wymagamy ponadto, aby po usunięciu ostatniego paska (a także dwóch ostatnich
pasków, trzech ostatnich pasków, itd.) pozostały kształt nadal był diagramem
Younga.

Dla każdego paska jego \emph{wysokość} zdefiniowana jest jako
różnica między
największym i najmniejszym numerem wiersza wszystkich jego klatek. W
przykładzie z Rysunku~\ref{rys:murnaghan} wysokości pasków wynoszą kolejno:
$1,0,1$. Wysokość rozkładu na paski definiujemy jako sumę wysokości wszystkich
pasków. Reguła Murnaghana-Nakayamy głosi, że wkład danego rozkładu na paski
wynosi $(-1)^{\text{(wysokość rozkładu)}}$. Aby obliczyć charakter $\Tr
\rho^\lambda(\pi)$ powinniśmy zsumować wkła\-dy od wszystkich rozkładów na
paski.

Prostym przykładem zastosowania reguły Murnaghana-Nakayamy jest analiza
nieredukowalnych reprezentacji grupy $\Sym{3}$.
Jak łatwo można się przekonać, istnieją dokładnie trzy diagramy Younga o trzech
klatkach; odpowiedniość pomiędzy nimi a reprezentacjami
nieredukowalnymi opisanymi w Rozdziale \ref{subsec:s3} może być wyznaczona
właśnie dzięki regule Murnaghana-Nakayamy poprzez porównanie charakterów.

Nieco trudniejszy przykład dotyczy reprezentacji grupy $\Sym{5}$. Reprezentację
$\rho$ grupy $\Alt{5}$ na przestrzeni $\R^3$ przedstawioną w Rozdziale
\ref{subsec:a5} można wykorzystać do konstrukcji interesującej reprezentacji
$\widetilde{\rho}$ grupy $\Sym{5}$ na przestrzeni $\R^6=\R^3\oplus \R^3$.
Mianowicie wybierzmy dowolną permutację nieparzystą $\sigma\in \Sym{5}$ o tej
własności, że $\sigma^2=e$. 
Żądamy, aby reprezentacja $\widetilde{\rho}$ miała dwie własności:
po pierwsze, aby dla dowolnej permutacji parzystej $\pi\in \Alt{5}$ zachodziło
$$ \widetilde{\rho}(\pi)(v_1,v_2) = \big( \rho(\pi)(v_1), \rho(\pi)(v_2)
\big);$$
innymi słowy obcięcie $\widetilde{\rho}$ do podgrupy $\Alt{5}$ ma być równe
$\rho\oplus\rho$.
Po drugie,
$$ \widetilde{\rho}(\sigma)(v_1,v_2) = (v_2,v_1). $$
Powyższe warunki okazują się jednoznacznie wyznaczać reprezentację
$\widetilde{\rho}$. Korzystając z charakterów obliczonych przy pomocy reguły
Murnaghana-Nakayamy można zidentyfikować ją z nieredukowalną reprezentacją
odpowiadającą diagramowi Younga $(3,1,1)$. Mam nadzieję, że powyższy przykład
przekona Czytelnika o tym, że reprezentacje grup permutacji mogą mieć bardzo
nietrywialną strukturę.

\subsection{Skalowanie i uogólnione diagramy Younga}
\label{sec:nastepny-po-murnaghan}

Pierwszy asymptotyczny problem, który chciałbym przedstawić, może być w 
nieformalny sposób sformułowany następująco: \emph{przypuśćmy, że
dany jest ciąg diagramów Younga $(\lambda^{(n)})$, który w jakimś sensie dąży do
nieskończoności; co możemy powiedzieć o odpowiadającym mu ciągu reprezentacji
nieredukowalnych $\rho^{\lambda^{(n)}}$?} Aby nadać temu problemowi sens,
spróbuję teraz doprecyzować w jakim sensie diagramy Younga mają dążyć do
nieskończoności. Można to zrobić na wiele sposobów, na razie skoncentruję się
na przypadku, w którym diagramy te, dążąc do nieskończoności, zachowują w
pewnym sensie swój kształt; poniżej przedyskutuję ten pomysł.

Jeśli $s>0$ jest liczbą całkowitą, zaś $\lambda$ jest diagramem Younga, przez
$s\lambda$ oznaczać będę diagram $\lambda$ przeskalowany o czynnik $s$. W
geometrycznym podejściu diagram $s\lambda$ powstaje przez podziałanie
jednokładnością o skali $s$ na diagram $\lambda$ lub, innymi słowy, przez
zastąpienie każdej z klatek diagramu $\lambda$ kratką złożoną z $s\times s$
klatek. Wynika stąd, że diagram $s\lambda$ składa się z $s^2 |\lambda|$ klatek.
W
niegeometrycznym podejściu
$$ s\lambda=(\underbrace{s\lambda_1,\dots,s\lambda_1}_{s \text{
razy}},\underbrace{s\lambda_2,\dots,s\lambda_2}_{s \text{ razy}},\dots )$$

Ponieważ interesuje mnie sytuacja, w której liczba klatek diagramu dąży do
nieskończoności, ale diagramy zachowują wspólny kształt, dobrym pomysłem jest
badanie ciągu diagramów $\lambda,2\lambda,3\lambda,\dots$ powstałego przez
przeskalowanie ustalonego diagramu $\lambda$.

Alternatywnym, nieco ogólniejszym podejściem, jest dopuszczenie skalowań
$s\lambda$, w których $s>0$ jest dowolną liczbą rzeczywistą. Działanie przez
jednokładność w skali $s$ na diagram $\lambda$ zwykle nie daje w wyniku diagramu
Younga, ale pewien geometryczny obiekt, który nazywać będę \emph{uogólnionym
diagramem Younga}. Korzystając z tego pojęcia można powiedzieć precyzyjnie co
oznacza, że ciąg diagramów Younga $\lambda^{(1)},\lambda^{(2)},\dots$ dąży do
nieskończoności dążąc do pewnego asymptotycznego kształtu: mianowicie liczba
klatek diagramu $\lambda^{(n)}$ ma dążyć do nieskończoności w granicy
$n\to\infty$ oraz ciąg uogólnionych diagramów Younga
$\frac{1}{\sqrt{|\lambda^{(n)}|}} \lambda^{(n)}$ ma zbiegać (w jakiejś rozsądnie
wybranej topologii) do jakiegoś uogólnionego diagramu Younga.

Jak widać, tego typu skalowanie odpowiada sytuacji, w której diagram Younga
$\lambda$ ma co najwyżej $C\sqrt{|\lambda|}$ wierszy i kolumn, gdzie liczba $C$
jest ustalona. Tego typy diagramy nazywane są $C$-zbalansowanymi. Możliwe jest
oczywiście badanie asymptotyki diagramów Younga w innych skalowaniach niż tylko
skalowanie zbalansowanych diagramów Younga; wykracza to jednak poza ramy
niniejszego artykułu. 

\subsection{Asymptotyka i algorytmy kombinatoryczne}

W dalszym ciągu używać będę konwencji, że dla $n\leq m$ każda permutacja
permutacja zbioru $\{1,2,\dots,n\}$ może być również traktowana jako
permutacja zbioru $\{1,2,\dots,n,\dots,m\}$; po prostu deklaruję, że ma ona
działać jako identyczność na dodatkowych elementach $\{n+1,\dots,m\}$.

Nareszcie mogę sformułować w konkretny sposób pierwszy problem asymptotycznej
teorii reprezentacji:
\begin{problem}
Załóżmy, że ciąg diagramów Younga $(\lambda^{(n)})$ w jakiś sposób dąży do
nieskończoności, a $\pi$ jest pewną ustaloną permutacją. Co możemy powiedzieć o
asymptotyce charakterów $\Tr \rho^{\lambda^{(n)}}(\pi)$ w granicy $n\to\infty$?
\end{problem}
Nie jest trudno wyobrazić sobie, że rozwiązanie powyższego problemu przy pomocy
reguły Murnaghana-Nakayamy jest zadaniem karkołomnym:  liczba sposobów podziału
diagramu Younga na paski rośnie błyskawicznie z rozmiarem diagramu, różne
podziały dają wkład z przeciwnymi znakami, ich wkłady często wzajemnie się
znoszą\dots Tego typu sytuacja jest dość typowa w asymptotycznej teorii
reprezentacji: na niemal każde pytanie znana jest od dawna odpowiedź w postaci
kombinatorycznego algorytmu, który jednak ze wzrostem rozmiaru problemu staje
się tak skomplikowany i nużący, że nie daje zbyt dokładnej informacji
na temat asymptotyki. Konieczne są więc inne, bardziej analityczne metody.

\subsection{Jak opisać kształt diagramu?}

Zwykły opis diagramu Younga jako słabo malejącego ciągu $\lambda=(\lambda_1,
\lambda_2,\dots)$ jest opisem zbyt dyskretnym i kombinatorycznym; z tego powodu
kiepsko nadaje się do stosowania w problemach asymptotycznych. 

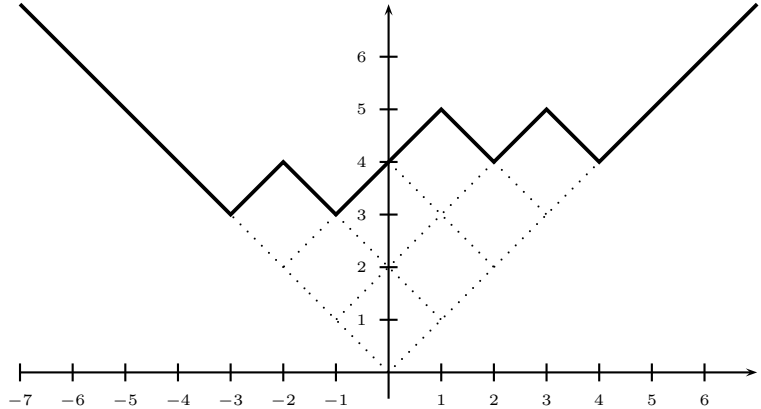
\begin{figure}
\begin{center}
\psset{unit=0.7cm}
\begin{pspicture}(-7.2,-0.2)(7.2,7.2)
\tiny
\psaxes{->}(0,0)(-7,-0.5)(7,7)
\psline[linewidth=0.5mm](7,7)(4,4)(3,5)(2,4)(1,5)(-1,3)(-2,4)(-3,3)(-7,7)
\psline[linestyle=dotted](-1,1)(2,4)
\psline[linestyle=dotted](-2,2)(-1,3)
\psline[linestyle=dotted](1,1)(-1,3)
\psline[linestyle=dotted](2,2)(0,4)
\psline[linestyle=dotted](3,3)(2,4)
\pspolygon[linestyle=dotted](0,0)(4,4)(3,5)(2,4)(1,5)(-1,3)(-2,4)(-3,3)
 
\end{pspicture}
\end{center}
 \caption{Diagramu Younga $(4,3,1)$ z Rys.~\ref{rys:francuski} narysowany w
konwencji
rosyjskiej. Gruba linia przedstawia \emph{profil} diagramu.}
\label{rys:rosyjska}
\end{figure}

Znacznie lepszym pomysłem jest użycie \emph{konwencji rosyjskiej} do rysowania
diagramów Younga. Diagram $(4,3,1)$ został narysowany w konwencji francuskiej
na Rysunku \ref{rys:francuski} oraz w konwencji rosyjskiej na
Rysunku \ref{rys:rosyjska}. Jak widać, rysunek w konwencji rosyjskiej powstaje z
rysunku w konwencji francuskiej przez złożenie obrotu o kąt $\frac{\pi}{4}$
oraz jednokładności o skali $\sqrt{2}$. Brzeg diagramu Younga, na
Rysunku~\ref{rys:rosyjska} zaznaczony pogrubioną linią, nazywany jest
\emph{profilem}
diagramu i może być utożsamiony z pewną funkcją $\omega:\R\rightarrow\R_+$.

Konwencja rosyjska ma wiele zalet. Przede wszystkim,  w świetle pracy Okounkova
i Vershika \cite{OkounkovVershik1996} jest bez wątpienia bardzo naturalną
konwencją rysowania diagramów Younga. Ponadto pozwala zdefiniować uogólnione
diagramy Younga jako ,,ciągłe profile'', czyli jako funkcje
$\omega:\R\rightarrow\R_+$ spełniające pewne proste aksjomaty.

Drugim dobrym pomysłem na opis kształtu diagramu Younga jest wprowadzenie
jakiejś nowej rodziny parametrów $(U_i)$, która zachowywałaby się w sposób mniej
kombinatoryczny, a bardziej analityczny. Możnaby wymienić tutaj wiele
pożądanych własności, jakich należy wymagać od takiej rodziny parametrów, aby
była ona jak najbardziej użyteczna; w niniejszym artykule wspomnę jedynie o
warunku \emph{jednorodności}. Mianowicie wymagać będę aby parametr $U_i$
był jednorodną funkcją diagramu Younga, stopnia $i$, to znaczy aby
$$ U_i(s\lambda) = s^i U_i(\lambda).$$
Jednorodność parametrów będzie ich wielkim atutem, gdyż pozwoli
bardzo dokładnie wyznaczać asymptotykę funkcji wielomianowych zależnych od
parametrów.

Jest kilka dobrych sposobów na wprowadzenie takiej rodziny parametrów i każdy
z nich wykazuje swoje zalety w nieco innej sytuacji. Aby dać Czytelnikowi
smak tego, jak taka rodzina może wyglądać, wspomnę o bardzo naturalnej
rodzinie parametrów $(S_i)$ zadanych wzorem \cite{DolegaF'eray'Sniady2008}
$$ S^\lambda_i = (i-1) \iint_{(x,y)\in\lambda} 
(x-y)^{i-2}\ dx\ dy  $$
dla diagramu Younga $\lambda$ narysowanego w układzie współrzędnych w konwencji
francuskiej i dla całkowitego $i\geq 2$. W przypadku, jeśli ten sam
diagram Younga zostanie narysowany w układzie współrzędnych w konwencji
rosyjskiej, parametry te wyrażają się wzorem
$$ S^\lambda_i = (i-1) \frac{1}{2} \iint_{(x,y)\in\lambda} 
x^{i-2}\ dx\ dy,  $$
innymi słowy
\begin{equation} 
\label{eq:S}
 S_i^\lambda = (i-1) \int_{-\infty}^{\infty} x^{i-2}\
\frac{\omega(x)-|x|}{2}
\ dx, 
\end{equation}
gdzie $\omega$ jest profilem diagramu Younga $\lambda$. Jak widać z ostatniej
równości, w konwencji rosyjskiej parametry te wyrażają się w wyjątkowo prosty
sposób, co potwierdza użyteczność konwencji rosyjskiej.

Parametr $S_2$ równy jest po
prostu liczbie klatek diagramu. Pozostałe parametry w mniejszym lub
większym stopniu również można zinterpretować geometrycznie: parametr
$S_3$ jest w pewnym sensie miarą tego jak bardzo nachylony jest profil
diagramu (w konwencji rosyjskiej) względem poziomu, a paramter $S_4$ opisuje
jak bardzo profil jest wygięty w kształt symbolu $\cup$, a jak bardzo w kształt
symbolu $\cap$.

\subsection{Znormalizowane charaktery}

Dla permutacji $\pi\in \Sym{k}$ oraz diagramu Younga $\lambda$ zadającego
reprezentację nieredukowalną grupy $\Sym{n}$ \emph{znormalizowany
charakter} zdefiniowany jest wzorem
$$ \Sigma_{\pi}^{\lambda} = \underbrace{n(n-1)\cdots(n-k+1)}_{k \text{
czynników}} \frac{\Tr \rho^{\lambda}(\pi)}{\Tr \rho^{\lambda}(e)}. $$ Iloraz
$\frac{\Tr \rho^{\lambda}(\pi)}{\Tr \rho^{\lambda}(e)}$ jest bardzo naturalną
wielkością, już raz napotkaliśmy go w równaniu \eqref{eq:normalizowane} podczas
badania spacerów losowych, z kolei iloczyn $n(n-1)\cdots(n-k+1)$ może być
interpretowany jako kombinatoryczna wielkość
opisująca na ile sposobów permutacja $\pi\in \Sym{k}$ może być zanurzona do
grupy
$\Sym{n}$. Ponieważ występujący w mianowniku charakter $\Tr \rho^{\lambda}(e)$
obliczony na elemencie neutralnym grupy jest równy po prostu wymiarowi
reprezentacji, może być on bezpośrednio wyliczony na podstawie
tzw.~\emph{formuły Robinsona-Thralla}. Podsumowując: znormalizowane charaktery
są naturalnymi wielkościami, które zawierają zasadniczo te same informacje, co
zwykłe charaktery. Jak zobaczymy dalej, znormalizowane charaktery szczególnie
dobrze nadają się do badania problemów asymptotycznych.

\subsection{Formuła Stanleya-F\'eraya}

Obecnie dostępnych jest kilka metod badania (znormalizowanych) charakterów grup
permutacji, które mniej lub bardziej nadają się do badania problemów
asymptotycznych. W niniejszym artykule przedstawię najmłodszą z nich, opartą na
\emph{formule Stanleya-F\'eraya}, gdyż ma ona rozliczne zalety: jest zarazem
najprostsza oraz, najprawdopodobniej, najpotężniejsza. Jest ona ponadto
przykładem nowych metod kombinatorycznych używanych przez asymptotyczną teorię
reprezentacji.

Przez $C(\pi)$ oznaczać będę zbiór cykli permutacji $\pi$. Dla danych
permutacji $\sigma_1,\sigma_2\in \Sym{k}$ będę rozważać kolorowanie $h$ 
cykli permutacji $\sigma_1$ (kolor każdego cyklu to numer pewnej kolumny
diagramu $\lambda$) oraz cyklów permutacji $\sigma_2$ (kolor każdego cyklu to
numer pewnego wiersza diagramu $\lambda$). Ściśle rzecz biorąc, takie
pokolorowanie jest funkcją na sumie rozłącznej zbiorów cykli permutacji
$\sigma_1$ i $\sigma_2$, a zatem 
$h:C(\sigma_1) \sqcup C(\sigma_2)\rightarrow\N$.
Mówimy, że \emph{kolorowanie $h$
jest zgodne z diagramem $\lambda$} jeśli dla wszystkich par cykli
$c_1\in C(\sigma_1)$ i $c_2\in C(\sigma_2)$ jeśli cykle $c_1$ i $c_2$ nie są
rozłączne, to na przecięciu kolumny 
$h(c_1)$ i wiersza $h(c_2)$ znajduje się klatka należąca do diagramu $\lambda$.
Przez $N^{\lambda}(\sigma_1,\sigma_2)$ będziemy oznaczać liczbę pokolorowań
cykli permutacji $\sigma_1$ i $\sigma_2$ zgodnych z diagramem $\lambda$. W
dalszej części rozdziału przedstawię przykład, który powinien
rozwiać wątpliwości związane  z powyższymi definicjami.

Następujące twierdzenie zostało sformułowane w nieco inny, równoważny
sposób, jako hipoteza przez Stanleya \cite{Stanley-preprint2006} i udowodnione
po raz pierwszy przez F\'eraya \cite{F'eray-preprint2006}, z tego powodu
nazywane jest ono formułą Stanleya-F\'eraya. Jego bardziej elementarny dowód 
można znaleźć w pracy \cite{F'eray'Sniady-preprint2007}.
\begin{theorem}
\label{theo:character-sniady}
Dla dowolnej permutacji $\pi\in \Sym{k}$ i dowolnego diagramu Younga $\lambda$
znormalizowany charakter $\Sigma^{\lambda}_{\pi}$ wyraża się następująco:
\begin{equation}
\label{eq:main-theorem}
 \Sigma^{\lambda}_{\pi}= \sum_{\substack{\sigma_1,\sigma_2\in \Sym{k},\\
\sigma_1
\sigma_2=\pi}}
(-1)^{\sigma_1}\ N^{\lambda}_{\sigma_1,\sigma_2},
\end{equation}
gdzie $(-1)^{\sigma_1}$ oznacza znak permutacji $\sigma_1$.
\end{theorem}

Jako przykład obliczę wartość znormalizowanego charakteru na
transpozycji $\pi=(12)$. Jeden z dwóch składników w sumie
\eqref{eq:main-theorem} odpowiada parze permutacji $\sigma_1=(1)(2)$ oraz
$\sigma_2=(12)$. Pokolorowania, które liczą się do
$N^{\lambda}_{\sigma_1,\sigma_2}$, są więc następującej postaci: jedynemu
cyklowi permutacji $\sigma_2=(12)$ odpowiada numer wiersza, zaś każdemu z dwóch
cykli permutacji $\sigma_1=(1)(2)$ odpowiada numer kolumny diagramu Younga.
Z warunku na zgodność wynika więc, że
$$ N^{\lambda}_{ (1)(2), (12) } = \sum_i (\lambda_i)^2, $$
gdzie $\lambda_i$ to liczba klatek w $i$-tym wierszu diagramu.
W podobny sposób można wykazać, że
$$ N^{\lambda}_{ (12), (1)(2) } = \sum_i (\lambda'_i)^2, $$
gdzie $\lambda'_i$ to liczba klatek w $i$-tej kolumnie diagramu.
Z Twierdzenia \ref{theo:character-sniady} wynika więc, że
$$\Sigma^{\lambda}_{(12)} = n(n-1) 
\frac{\Tr \rho^{\lambda}_{(12)}}{\Tr \rho^{\lambda}(e)}= 
N^{\lambda}_{ (1)(2), (12) }-N^{\lambda}_{ (12), (1)(2)}=
\sum_i (\lambda_i)^2-\sum_i (\lambda'_i)^2.$$

Wzór Stanleya-F\'eraya jest interesujący z kilku powodów. Po
pierwsze, jego stopień komplikacji  (czyli liczba składników z prawej
strony) zależy tylko od permutacji $\pi$ na której obliczamy znormalizowany
charakter, a nie zależy od rozmiaru diagramu Younga $\lambda$. Z tego powodu
doskonale nadaje się do badania asymptotyki, w której permutacja $\pi$ jest
ustalona, zaś diagram Younga $\lambda$ dąży do nieskończoności. Po drugie, jest
to formuła o charakterze raczej kombinatorycznym, a nie analitycznym, gdyż suma
przebiega po \emph{faktoryzacjach permutacji $\pi$}, to znaczy po rozwiązaniach
równania $\sigma_1 \sigma_2=\pi$, a badanie różnego rodzaju faktoryzacji
permutacji jest jedną z ważnych dziedzin kombinatoryki.

Nie sposób nie wspomnieć o pewnych wadach wzoru Stanleya-F\'eraya---przede
wszystkim jest on raczej narzędziem teoretycznym niż praktycznym, gdyż liczba
składników w nim występujących szybko rośnie ze wzrostem komplikacji
permutacji $\pi\in\Sym{k}$ i z tego powodu do implementacji komputerowej należy
używać innych metod.


\subsection{Asymptotyka znormalizowanych charakterów}
\label{sec:poptokach}

Dla ustalonych permutacji $\sigma_1,\sigma_2\in\Sym{k}$ porównajmy liczbę
ich pokolorowań zgodnych z diagramem $\lambda$ oraz liczbę pokolorowań zgodnych
z przeskalowanym diagramem $s\lambda$ dla
całkowitego $s\geq 1$. Jeśli $h$ jest jakimś pokolorowaniem zgodnym z
przeskalowanym diagramem
$s\lambda$, to kolorowanie $F(h)$ zadane wzorem
$$ \big( F(h) \big)(c) = \left\lceil \frac{h(c)}{s} \right\rceil, $$
gdzie  $\left\lceil x \right\rceil$ oznacza zaokrąglenie w górę do najbliższej
liczby całkowitej, jest pokolorowaniem zgodnym z oryginalnym diagramem
$\lambda$. Ponieważ
każde pokolorowanie zgodne z $\lambda$ jest obrazem dokładnie
$s^{|C(\sigma_1)|+|C(\sigma_2)|}$ pokolorowań zgodnych z $s\lambda$, gdzie
$|C(\sigma_1)|+|C(\sigma_2)|$ oznacza łączną liczbę cykli permutacji $\sigma_1$
i $\sigma_2$, zatem
$$ N^{s\lambda}_{\sigma_1,\sigma_2} = s^{|C(\sigma_1)|+|C(\sigma_2)|}
N^{\lambda}_{\sigma_1,\sigma_2}. $$
Fakt ten możemy sformułować następująco: liczba pokolorowań jest
\emph{jednorodną} funkcją diagramu Younga, stopnia
$|C(\sigma_1)|+|C(\sigma_2)|$.

Z Twierdzenia \ref{theo:character-sniady} płynie więc szczególnie
prosty dowód następującego zaskakującego wniosku (który znany był na długo
przed Twierdzeniem \ref{theo:character-sniady}):
\begin{wniosek}
\label{wniosek:wielomian}
Dla dowolnej permutacji $\pi$ oraz dla dowolnego diagramu Younga $\lambda$
funkcja
$$ s\mapsto \Sigma_{\pi}^{s\lambda}$$ 
jest wielomianem.
\end{wniosek}

Nasz cel, zbadanie asymptotyki charakterów, sprowadza się więc do zbadania
współczynników powyższych wielomianów.


\subsection{Wolne kumulanty}

Używać będę specjalnej notacji do oznaczania znormalizowanych charakterów
na cyklach, mianowicie definiuję $$ \Sigma_k^\lambda:=
\Sigma^\lambda_{(1,2,\dots,k)},$$
gdzie $(1,2,\dots,k)$ traktuję jako element grupy $\Sym{k}$.

Korzystając z formuły Stanleya-F\'eraya nie jest trudno wykazać, że wielomian 
$$ s\mapsto \Sigma_{k-1}^{s\lambda}$$ 
jest wielomianem stopnia $k$. Dla $k\geq 2$ całkowitego
definiujemy \emph{$k$-tą wolną kumulantę} diagramu $\lambda$ jako
współczynnik stojący przy najwyższej potędze $s$:
\begin{equation}
\label{eq:free}
R^{\lambda}_k = [s^k] \Sigma_{k-1}^{s\lambda}. 
\end{equation}

Wolne kumulanty mają rozliczne zalety. Po pierwsze są jednorodne (to znaczy:
$R^{s\lambda}_k= s^k R^{\lambda}_k$), a zatem diagramy Younga o różnych
rozmiarach a tym samym kształcie mają te same wolne kumulanty (z dokładnością
do prostego przeskalowania), co bardzo ułatwia badanie asymptotyki funkcji
wyrażonych jako wielomiany w wolnych kumulantach. Po drugie, z samej definicji
\eqref{eq:free} wynika,
że wolne kumulanty dają asymptotykę charakterów na cyklach, co w nieco
nieformalny sposób można
zapisać jako
\begin{equation}
\label{eq:approx}
\Sigma_k^{\lambda} \approx R^{\lambda}_{k+1}.
\end{equation}
Wreszcie wolne kumulanty \emph{dają się efektywnie wyliczyć}. Z powyższych
powodów wolne kumulanty należą do najbardziej popularnych parametrów
opisujących kształt diagramu Younga. 

Ostatnia z wymienionych własność efektywnej wyliczalności jest
szczególnie ważna, gdyż jaki byłby pożytek z wielkości, której nie dałoby się
w praktyce obliczyć dla konkretnego diagramu Younga? Jest wiele
sposobów obliczania wolnych kumulant, w niniejszym artykule wspomnę tylko o
dwóch. Jak się okazuje, wolne kumulanty wyrażają się przez bardzo proste
parametry $(S_i)$ zdefiniowane równaniem \eqref{eq:S} jak następuje:
$$ R_{n}^\lambda =
 \sum_{l\geq 1} \frac{1}{l!} (-n+1)^{l-1} \sum_{\substack{k_1,\dots,k_l\geq 2 \\
k_1+\cdots+k_l=n }} S^\lambda_{k_1} \cdots S^\lambda_{k_l}. $$

Druga metoda sprowadza się do obliczenia charakteru $\Sigma_{k-1}$ dzięki
formule Stanleya-F\'eraya i sprawdzeniu dla których par permutacji
$\sigma_1,\sigma_2$ liczba pokolorowań $N^{s\lambda}_{\sigma_1,\sigma_2}$ jest
jednomianem parametru $s$ stopnia $k$. Pary permutacji o tej własności
nazywane są \emph{minimalnymi faktoryzacjami cyklu} i okazują
się mieć bardzo piękną kombinatoryczną strukturę związaną z
tzw.~\emph{nieprzecinającymi się partycjami}. 

\subsection{Wielomiany Kerova}

To bardzo zaskakujące, ale wolne kumulanty nie tylko zadają przybliżone,
asymptotyczne wartości charakterów jak w formule \eqref{eq:approx}, ale mogą
one również być użyte do obliczania \emph{dokładnych wartości charakterów}.

Sergey Kerov podczas wykładu w Institut Henri Poincaré w Paryżu w 2000 roku
naszkicował dowód następującego twierdzenia:
\begin{theorem}
Dla każdej permutacji $\pi$ istnieje wielomian o współczynnikach
całkowitych $K_{\pi}$ (zwany obecnie \emph{wielomianem Kerova}) o tej własności,
że
$$ \Sigma_{\pi}^\lambda= K_{\pi}(R_2^\lambda,R_3^\lambda,\dots) $$
zachodzi dla dowolnego diagramu Younga $\lambda$. 
\end{theorem}
Wielomiany Kerova nazywane są wielomianami uniwersalnymi, ponieważ nie zależą
one od wyboru diagramu Younga $\lambda$.

W przypadku, gdy permutacja jest cyklem, używamy specjalnej notacji
$$ \Sigma_{k}^\lambda= K_{k}(R_2^\lambda,R_3^\lambda,\dots). $$
Dla uproszczenia pomijać będę zależność charakterów i wolnych
kumulant od diagramu Younga i zamiast powyższej równości pisać będę
$$ \Sigma_{k}= K_{k}(R_2,R_3,\dots). $$

Kilka pierwszych wielomianów Kerova zawiera poniższa tabela.
\begin{align*}
\Sigma_1 &= R_2, \\
\Sigma_2 &= R_3, \\
\Sigma_3 &= R_4 + R_2,   \\
\Sigma_4 &= R_5 + 5R_3,    \\     
\Sigma_5 &= R_6 + 15R_4 + 5R_2^2 + 8R_2, \\
\Sigma_6 &= R_7 + 35R_5 + 35R_3 R_2 + 84R_3.
\end{align*}
Na podstawie powyższego przykładu niemal nie sposób nie sformułować
następującego przypuszczenia, które postawił jako pierwszy Kerov podczas swego
wykładu.
\begin{conjecture}
Współczynniki wielomianów Kerova $K_k$ są \emph{nieujemnymi} liczbami
całkowitymi. 
\end{conjecture}

Przypuszczenie to zostało udowodnione niedawno przez F\'eraya
\cite{F'eray-preprint2008}, ponadto w pracy 
\cite{DolegaF'eray'Sniady2008} podano
interpretację kombinatoryczą współczynników wielomianów Kerova $K_k$ jako liczby
faktoryzacji $\sigma_1 \sigma_2=(1,2,\dots,k)$ o pewnych dodatkowych
własnościach, o których zainteresowany Czytelnik może przeczytać w przystępnie
napisanym wstępie do pracy
\cite{DolegaF'eray'Sniady2008}. Powyższy wynik można wyrazić następująco: udało
się uzyskać nową, kombinatoryczną formułę na znormalizowane charaktery
$\Sigma_\pi^\lambda$ wyrażoną w języku wolnych kumulant diagramu $\lambda$.
Jej przewaga nad wzorem Stanleya-F\'eraya staje się widoczna, gdy
zauważymy, że w tym ostatnim składniki występują z przeciwnymi znakami, w
związku z czym często ich wkłady kasują się wzajemnie, podczas gdy wielomian
Kerova oraz formuła na jego współczynniki zawierają znacznie mniej składników i
dają się lepiej kontrolować.

Wyniki zawarte w obu wyżej cytowanych pracach oparte są na analizie składników
występujących we wzorze Stanleya-F\'eraya, a zatem można o nich myśleć
jak o wzorze Stanleya-F\'eraya po uwzględnieniu możliwie dużej liczby skracań.
Mam nadzieję, że te nowe wyniki dotyczące wielomianów Kerova pozwolą w
przyszłości na lepsze zrozumienie struktury reprezentacji grup permutacji, a w
szczególności że pozwolą lepiej szacować wartości charakterów.

\subsection{Dalsza lektura}

Dobrym wstępem do teorii reprezentacji grup permutacji oraz kombinatorycznych
aspektów diagramów Younga jest książka Sagana \cite{SaganSymmetric}. Prostym
wprowadzeniem do wielomianów Kerova jest praca Biane'a \cite{BianeCharacters}.
Wstęp do pracy \cite{DolegaF'eray'Sniady2008} zawiera dość szerokie tło
historyczne związane z wielomianami Kerova.
O związkach minimalnych faktoryzacji cyklu i
nieprzecinających się partycji traktuje praca Biane'a \cite{Biane1997crossings}.

\section{Podsumowanie}

Asymptotycznej teorii reprezentacji daleko jest jeszcze do dojrzałości.
Na przykład, znanych jest co najmniej kilka wzorów na charaktery; są one
nierzadko sformułowane przy pomocy dość odległych od siebie dziedzin
matematyki, jak choćby analiza i kombinatoryka. To, że wyrażają one tę samą
wielkość wcale nie jest oczywiste. Co więcej, pewne własności symetrii
związanych z charakterami wielomianów Kerova stają się oczywiste przy użyciu
jednych wzorów, a inne przy zastosowaniu innych. Nie widać obecnie jakiegoś
sposobu na zunifikowanie tych wszystkich podejść, które pozwalałoby badać
wszystkie symetrie wielomianów Kerova naraz. 

W poprzednim rozdziale wspomniałem o rozwiązaniu hipotezy Kerova.
Ten postęp w zrozumieniu wielomianów Kerova oczywiście
bardzo cieszy, ale sformułowanych zostało kilka nowych hipotez dotyczących
struktury wielomianów Kerova oraz dodatniości ich współczynników, jeśli zamiast
wolnych kumulant użyć pewnych innych naturalnych wielkości (hipotezy te zostały
zebrane we wstępie do pracy \cite{DolegaF'eray'Sniady2008}). Pozostają one nadal
otwarte.

Mam nadzieję, że dzięki naszkicowanym w niniejszej pracy nowym metodom
kombinatorycznym w przyszłości uda się osiągnąć lepsze zrozumienie teorii
reprezentacji oraz jej związków z innymi działami matematyki.

\section{Podziękowania}

Praca naukowa finansowana ze środków na naukę w latach 2006--2009 jako projekt
badawczy Ministerstwa Nauki i Szkolnictwa Wyższego numer
P03A 013 30, a także przez EC Marie Curie Host Fellowship for the
Transfer of Knowledge ``Harmonic
Analysis, Nonlinear Analysis and Probability'', numer kontraktu
MTKD-CT-2004-013389.

%

\noindent
Piotr Śniady \\ 
Instytut Matematyczny, Uniwersytet Wrocławski, pl.~Grunwaldzki 2/4, \mbox{50-384
Wrocław} \\
e-mail: \texttt{piotr.sniady@math.uni.wroc.pl}

\end{document}